\documentclass[number,review,12pt]{elsarticle}

 \usepackage{graphicx}
 \usepackage{eurosym}
 \usepackage{bbm}
 \usepackage{setspace}
 \usepackage{booktabs}
 \usepackage{url}
 \usepackage[top=1.3in, bottom=1.5in, left=1.4in, right=1.4in]{geometry}

 \usepackage{amssymb}

\journal{Applied Energy}

\begin{document}

\begin{frontmatter}

\title{Optimal Charging of an Electric Vehicle using a Markov Decision Process}

\author[address1]{Emil B. Iversen \corref{cor1}}
\ead{jebi@dtu.dk}

\author[address1]{Juan M. Morales}
\ead{jmmgo@dtu.dk}

\author[address1]{Henrik Madsen}
\ead{hm@dtu.dk}

\address[address1]{DTU Compute, Technical University of Denmark, Matematiktorvet, building 322, DK-2800 Lyngby, Denmark}
\cortext[cor1]{Corresponding author. Tel.: +45 45 25 30 75}

\begin{abstract}
The combination of electric vehicles (EVs) and renewable energy is taking shape as a potential driver for a future free of fossil fuels. However, the efficient management of the EV fleet is not exempt from challenges. It calls for the involvement of all actors directly or indirectly related to the energy and transportation sectors, ranging from governments, automakers and transmission system operators, to the ultimate beneficiary of the change: the end-user. An EV is primarily to be used to satisfy driving needs, and accordingly charging policies must be designed primarily for this purpose. The charging models presented in the technical literature, however, overlook the stochastic nature of driving patterns. Here we introduce an efficient stochastic dynamic programming model to optimally charge an EV while accounting for the uncertainty inherent to its use. With this aim in mind, driving patterns are described by an inhomogeneous Markov model that is fitted using data collected from the utilization of an EV. We show that the randomness intrinsic to driving needs has a substantial impact on the charging strategy to be implemented.
\end{abstract}

\begin{keyword}
    Electric vehicles \sep
    Driving patterns \sep
    Optimal charging \sep
    Markov processes \sep
    Stochastic dynamic programming \sep
\end{keyword}

\end{frontmatter}





%


\section{Introduction}\label{intro}

Electric vehicles (EVs) are emerging as a sustainable and environmentally friendly alternative to conventional vehicles, provided that the energy used for their charging is obtained from renewable energy sources. The energy generated from renewable sources such as sunlight, wind and waves is, however, dependent on weather conditions. As a consequence, the electricity production from these sources is inherently uncertain in time and quantity. Furthermore, electricity has to be produced and consumed at the same time, as the large-scale storage of the energy generated is, still today, very limited. As a result, the energy obtained from renewables may be wasted in times when the demand for electricity is not high enough to absorb it, with a consequent detrimental effect on the profitability of renewables. Since the battery in an EV is basically a storage device for energy, the large-scale integration of EVs in the transportation sector may contribute to substantially increasing the socioeconomic value of an energy system with a large renewable component, while reducing the dependence of the transportation sector on liquid fossil fuel.

For this reason, EVs have received increased interest from the scientific community in recent years (detailed literature reviews of the state of the art can be found in \cite{Green11} and \cite{bessa2012economic}). Special attention has been given to the analysis of the effect of EVs integration on the electricity demand profile~\cite{Weiller20113766, Shao2012543}, emissions~\cite{sioshansi2012modeling} and social welfare \cite{Juul, shortt2009impact, kiviluoma2010influence}, and to the design of charging schemes that avoid increasing the peak consumption \cite{Gan11, Ma10}, help mitigate voltage fluctuations and overload of network components in distribution grids~\cite{Richardson12}, and/or get the maximum economic benefit from the storage capability of EVs within a market environment, either from the perspective of a single vehicle~\cite{Kempton2005268, Rotering11} or the viewpoint of an aggregator of EVs~\cite{kristoffersen2011optimal,Momber}. In all these publications, though, and more generally in the technical literature on the topic, the charging problem of an EV is addressed either by considering deterministic driving patterns, when the focus is placed on the management of a single vehicle, or by aggregating the driving needs of different EV users, when the emphasis is on modeling a whole fleet of EVs. This aggregation, however, obscures the dynamics of each specific vehicle. Likewise, the deterministic driving patterns of a single EV are often based on expected values or stylized behaviors, which fail to capture important features of the charging problem such as the daily variation in the use of the vehicle or potential user conflicts in terms of not having the vehicle charged and ready for use. A stochastic model for driving patterns provides more insight into these aspects and becomes fundamental for applying a charging scheme in the real world. Despite this, the stochastic modeling of driving patterns has received little attention from the scientific community, as pointed out in \cite{Green11}. We mention here the research work by \cite{Lojowska11}, in which they aim to capture the uncertainty intrinsic to the vehicle use by means of a Monte Carlo simulation approach. They assume, however, an uncontrolled charging scheme.

The work developed in this paper departs from the following two premises:

\begin{enumerate}
  \item The primary purpose of the battery of an EV is to provide power to drive the vehicle and not to store energy from the electricity grid. Consequently, it is essential that enough energy is kept in the battery to cover any desired trip. This calls for a decision tool that takes into account the driving needs of the EV user to determine when charging can be postponed and when the battery should be charged right away.
  \item The complexity of human behavior points to a stochastic model for describing the use of the vehicle. In turn, this stochastic model should be integrated into the aforementioned decision tool and exploited by it.
\end{enumerate}

That being so, this paper introduces an algorithm to optimally decide when to charge an EV that exhibits a stochastic driving pattern. The algorithm builds on the inhomogeneous Markov model proposed in \cite{Iversen12} for describing the stochastic use of a single vehicle. The model parameters are then estimated on the basis of data from the use of the specific vehicle. The approach captures the diurnal variation of the driving pattern and does not rely on any assumptions on the use of the vehicle, which makes it general and particularly versatile. Our algorithm thus embodies a \emph{Markov decision process} which is solved recursively using a stochastic dynamic programming approach. The resulting decision-support tool allows for addressing issues related to charging, vehicle-to-grid (V2G) schemes \cite{kempton1997electric, Kempton2005268}, availability and costs of using the vehicle. The algorithm runs swiftly on a personal computer, which makes it feasible to implement on an actual EV.

The remainder of this paper is organized as follows: In Section 2 the stochastic model for driving patterns developed in \cite{Iversen12} is briefly described, tailored to be used in the present work, and extended to address the problem of driving data limitations through hidden Markov models. Section 3 introduces the algorithm for the optimal charging of an EV as a Markov decision process that is solved using stochastic dynamic programming. Section 4 provides results from a realistic case study and explores the potential benefit of implementing V2G schemes. Section 5 concludes and provides directions for future research within this topic.

\section{A Stochastic Model for Driving Patterns} \label{Stochastic Model}

In this section we summarize and extend the stochastic model for driving patterns developed in \cite{Iversen12}. We refer the interested reader to this work for a detailed description of the modeling approach.

\subsection{Standard Markov Model}

A state-space model is considered to describe the use of the EV. In its simplest form, it contains two states, according to which the vehicle is either \emph{driving} or \emph{not driving}. A more extensive version of the model would include a larger number of states which could capture information about where the vehicle is parked, how fast it is driving or what type of trip it is on. The basics of the general multi-state stochastic model are described in this section, including how to fit a specific model on an observed data set.

Let $X_t$, where $t \in \{ 0,1,2, \ldots \}$, be a sequence of random variables that takes on values in the countable set $S$, called the state space. Denote this sequence as $X$. We assume a finite number, $N$, of states in the state space. A Markov chain is a random process where future states, conditioned on the present state, do not depend on the past states \citep{Grimmett01}. In discrete time $X$ is a Markov chain if
\begin{equation}
\mathbb{P}\left( X_{t+1} = k | X_0 =x_0 , \ldots , X_{t} =x_{t} \right) = \mathbb{P}\left( X_{t+1} = k | X_{t} =x_{t} \right)
\end{equation}
for all $t \geq 0$ and all $\{ k, x_0 , \ldots , x_{t} \} \in S$.

A Markov chain is uniquely characterized by the transition probabilities, $p_{jk}(t)$, i.e.
\begin{equation}
\label{transprob}    p_{jk}(t) = \mathbb{P} \left( X_{t+1} = k | X_t = j \right).
\end{equation}
If the transition probabilities do not depend on $t$, the process is called a homogeneous Markov chain. If the transition probabilities depend on $t$, the process is known as an inhomogeneous Markov chain.

 When it comes to the use of a vehicle, it is appropriate to assume that the probability of a transition from state $j$ to state $k$ is similar on specific days of the week. Thus, for instance, Thursdays in different weeks will have the same transition probabilities. For convenience we further assume that all weekdays (Monday through Friday) have the same transition probabilities. These assumptions can be easily relaxed or interchanged with other assumptions and as such, are not essential to the model. With a sampling time in minutes, and taking into account that there are 1440 minutes in a day, this leads to the assumption:
\begin{equation}
    p_{jk}(t) = p_{jk}(t+1440).   \label{periodicity}
\end{equation}
This assumption implies that the transition probabilities, defined by (\ref{transprob}), are constrained to be a function of the time, $s$, in the diurnal cycle. Let the matrix containing the transition probabilities be denoted by $\mathbf{P}(s)$. For the model containing $N$ states the transition probability matrix is given by:
 \begin{equation}
    \mathbf{P}(s) =\left(
            \begin{array}{cccc}
              p_{11}(s) & p_{12}(s) & \ldots & p_{1N}(s) \\
              p_{21}(s) & p_{22}(s) & \ldots & p_{2N}(s) \\
              \vdots &  \vdots & \ddots & \vdots \\
              p_{N1}(s) & p_{N2}(s) & \ldots & p_{NN}(s) \\
            \end{array}
          \right),
\end{equation}
where $p_{jj}(s) = 1 - \sum_{i=1,i \neq j}^N{p_{ji}}$.

Now let $n_{jk}(s)$ define the number of observed transitions from state $j$ to state $k$ at time $s$. From the conditional likelihood function, the maximum-likelihood estimate of $p_{jk}(s)$ can then be found as:
\begin{equation}
    \widehat{p}_{jk}(s) = \frac{n_{jk}(s)}{\sum_{k=1}^N{n_{jk}(s)}}. \label{P_hat}
\end{equation}

A discrete time Markov model can be formulated based on the estimates of $\mathbf{P}(1), \mathbf{P}(2), \ldots, \mathbf{P}(1440)$. One apparent disadvantage of such a discrete time model is its huge number of parameters, namely $N \times (N-1) \times 1440$, where $N \times(N-1)$ parameters have to be estimated for each time step. Needless to say, the number of parameters to be estimated increases as the number of states grows. We refer to \cite{Iversen12} for further details on techniques to reduce the number of parameters to be estimated for each time step for models with more than two states. Another problem is linked to the number of observations available to properly carry out the estimation, i.e. if $\sum_{k=1}^N{n_{jk}(s^{\prime})} = 0$ for some $s^{\prime}$, then $\widehat{p}_{jk}(s^{\prime})$ is undefined.

To deal with the large number of parameters as well as undefined transition probability estimates, B-splines are applied to capture the diurnal variation in the driving pattern through a \emph{generalized linear model}. The procedure of applying a \emph{generalized linear model} is implemented in the statistical software package R as the function \texttt{glm($\cdot$)}. For a thorough introduction to B-splines see \cite{Hastie08} and for a general treatment of generalized linear models see \cite{Madsen10}. Next we elaborate on how the fitting of the Markov chain model works in our particular case.

Each day, at a specific minute, a transition from state $j$ to state $k$ either occurs or does not occur. Thus for every $s$ on the diurnal cycle we can consider the number of transitions to be binomially distributed, i.e. $n_{jk}(s) \sim B(z_{j}(s),p_{jk}(s))$, where the number of Bernoulli trials at $s$, given by $z_{j}(s)=\sum_{k=1}^N{n_{jk}(s)}$, is known and the probability of success, $p_{jk}(s)$, is unknown. The data can now be analyzed using a \emph{logistic regression}, which is a generalized linear model \citep{Madsen10}. The explanatory variables in this model are taken to be the basis functions for the B-spline. The logit transformation of the odds of the unknown binomial probabilities are modeled as linear combinations of the basis functions. We model $Y_{jk}(s) = n_{jk}(s) / z_{j}(s)$ and in particular, we are interested in $\textrm{E}\left[ Y_{jk}(s) \right] = p_{jk}(s)$.

As the basis functions for the B-spline are uniquely determined by the knot vector $\tau$, deciding the knot position and the amount of knots is important to obtain a good fit for the model. Here we proceed as follows: First a number of knots are placed on the interval $\left[ 0, 1440 \right]$, with one at each endpoint and equal spacing between them. Denote this initial vector of knots by $\tau_{\textrm{ \scriptsize init}}$. The model is then fitted using the basis functions as explanatory variables. Next, the fit of the model between the knots is evaluated via the likelihood function and an additional knot is placed in the center of the interval with the lowest likelihood value. The new knot vector is then given by $\tau'$. We repeat this procedure until the desired number of knots is reached. To determine the appropriate number of knots and avoid over-parametrization, on the basis of a likelihood ratio principle, we test that adding a new knot does significantly improve the fit.

\subsection{Hidden Markov Models}

Standard Markov models are limited in the sense that only those states that are actually observed can be modeled. Thus, if the data at our disposal only provide information on when the vehicle is either \emph{driving} or \emph{not driving}, the standard Markov model is restricted to two states. Furthermore the time spent in each state is exponentially distributed, albeit with time-varying intensity, and accordingly, the time until the next transition does not depend on the time spent in the current state. This may be particularly unrealistic for a model with few states for describing driving patterns.

To overcome these limitations, we can use a hidden Markov model, which allows estimation of additional states that are not directly observed in the data. In fact, we can estimate these states so that the waiting time in each state matches that which is actually observed in the data. Adding a hidden state is done by introducing a new state in the underlying Markov chain. The new state, however, is indistinguishable from any of the previously observed states. This allows for the waiting time in each observable state to be the sum of exponential variables, which is a more versatile class of distributions. It is worth insisting that the use of hidden Markov models is justified here to address insufficient state information in our data, which only include whether the vehicle is \emph{driving} or \emph{not driving}. Indeed, the same results could be obtained using the underlying Markov chain without hidden states, provided that the hidden states could be observed. In practice, though, more detailed driving data (e.g. including driving speed and/or location of the vehicle) could be available once the actual implementation is made on a vehicle, which in turn would avert the need for a hidden Markov model. For a detailed introduction to hidden Markov models, see \cite{zucchini2009hidden}, where techniques and scripts for estimating parameters are also provided.

The hidden Markov model consists of two parts. Firstly, an underlying unobserved Markov process, $\left\{X_t : t = 1,2, \ldots \right\}$, which describes the actual state of the vehicle. This part corresponds to the Markov model with no hidden states as described previously. The second part of the model is a state-dependent process, $\left\{ Z_t : t = 1,2, \ldots \right\}$, such that when $X_t$ is known, the distribution of $Z_t$ depends only on the current state $X_t$. A hidden Markov model is thus defined by the state-dependent transition probabilities, $p_{jk}(t)$, as defined for the standard Markov chain and the state-dependent distributions given by (in the discrete case):
\begin{equation}
\label{statedist}    d_{zk}(t) = \mathbb{P} \left( Z_{t} = z| X_t = k \right).
\end{equation}
Collecting the $d_{zk}(t)$'s in the matrix $\mathbf{D}(z_t)$, the likelihood of the hidden Markov model is given by:
\begin{equation}
\label{eq: HMM likelihood}    L_T = \delta \mathbf{D}(z_1) P(2) \mathbf{D}(z_2) \ldots P(T) \mathbf{D}(z_T),
\end{equation}
where $\delta$ is the initial distribution of $X_1$. We can now maximize the likelihood of observations to find the estimates of the transition probabilities.

\subsection{Fitting the Data}

The data at our disposal is from the utilization of a single vehicle in Denmark in the period spanning the six months from 23-10-2002 to 24-04-2003, with a total of 183 days. The data is GPS-based and follows specific cars. One car has been chosen and the model is intended to describe the use of this vehicle accordingly. The data set only contains information on whether the vehicle was \emph{driving} or \emph{not driving} at any given time. No other information was provided in order to protect the privacy of the vehicle owner. The data is divided into two periods, a training period for fitting the model from 23-10-2002 to 23-01-2003, and a test period from 24-01-2003 to 24-04-2003 for evaluating the performance of the model. The data set consists of a total of 749 trips. The time resolution is in minutes.

We shall consider a model with one \emph{not driving} state and several (hidden) \emph{driving} states. That is, it is observable if the vehicle is driving, but not the specific driving state. To fit the model to the data, we assume that only the transition probability from the \emph{not driving} state depends on the time of day. This is done to reduce the complexity of the estimation procedure, as it is cumbersome to estimate the time-varying parameters of a hidden Markov model. It is worth noting that a hidden Markov model allows for the probability of ending the current trip to depend on the time since departure, as the vehicle may pass through different driving states before ending the trip.

We now elaborate on the fitting of the hidden Markov model, which is split into estimation of its time-varying and time-invariant parameters.

\subsubsection{Fitting Time-Varying Parameters}
We need to estimate the probability of a transition from the vehicle being parked to a driving state. We denote this transition estimate by $\widehat{p}_{1 \cdot}(s)$. It holds that $\widehat{p}_{1 \cdot}(s) = 1 - \widehat{p}_{11}(s)$. Since both the parked state and the transitions from it are directly observable in the data, we can use the procedure described in Section 2.1 to estimate $\widehat{p}_{1 \cdot}(s)$.

\begin{figure}[h!]
\begin{center}
  \includegraphics[width=0.8\textwidth]{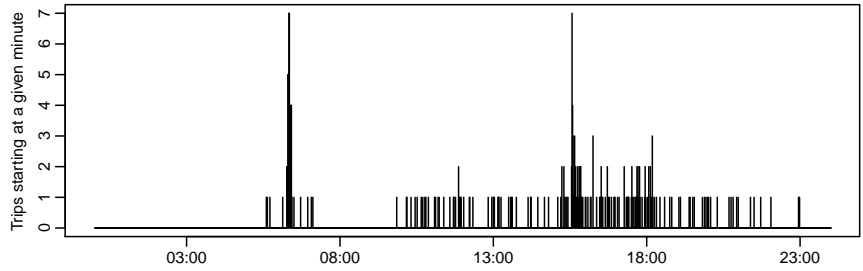}\\
    \vspace{-10pt}
  \caption{\emph{Number of trips starting at a certain minute of the day, cumulated for the first 66 weekdays.}}\label{StartingWeekdays}
\end{center}
\vspace{-10pt}
\end{figure}

The data have been divided into two main periods: weekdays and weekends. The observed number of trips starting every minute for the weekdays is displayed in Fig. \ref{StartingWeekdays}. A high degree of diurnal variation is found, with a lot of trips starting around 06:00 and again around 16:00. Also, there are no observations of trips starting between 00:00 and 05:00. Other patterns are found for weekends, but as these do not involve any methodological difference, we limit ourselves to trips starting on weekdays. Annual variations may also be present, however the limited data sample does not allow for capturing such seasonality.

\begin{figure}[h!]
  \begin{center}
    \includegraphics[width=0.8\textwidth]{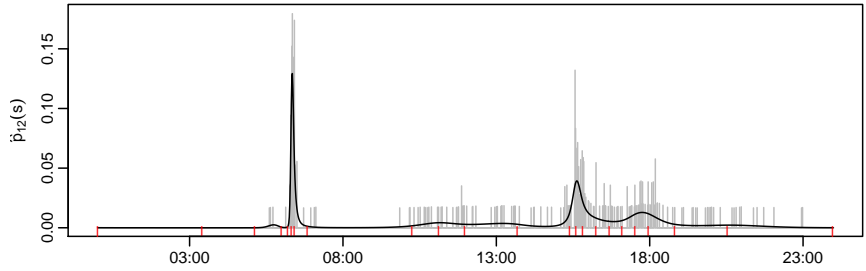}
  \end{center}
  \vspace{-20pt}
  \caption{\emph{$\widehat{p}_{1 \cdot}(s)$ based on the B-splines and the logistic regression, plotted as the black line over the estimates $\widehat{p}_{1 \cdot}(s)$ from (\ref{P_hat}), in gray. The red bars indicate the knot positioning.}}
  \label{TransitionProbability}
  \vspace{-10pt}
\end{figure}

The plot in Fig. \ref{TransitionProbability} illustrates the estimate of $\widehat{p}_{1 \cdot}(s)$ using B-splines with eight initial knots placed uniformly on the interval and 22 knots in total.

\subsubsection{Fitting Time-Invariant Parameters}

The time-invariant parameters are to be estimated so that an appropriate probability distribution is fitted to the duration of the trips. The time-invariant parameters are estimated by maximizing the likelihood given in (\ref{eq: HMM likelihood}). For a given number of \emph{driving} states, the transition probabilities can be estimated using the approach in \cite{zucchini2009hidden}. Once a model with $N$ states is fitted, we can test if adding an additional state significantly improves the fit. As a model with $N$ states is a sub-model of one with $N+1$ or more states, we increase the number of states until no significant improvement test is observed according to the likelihood ratio.

\begin{figure}[h!]
  \begin{center}
    \includegraphics[width=0.5\textwidth]{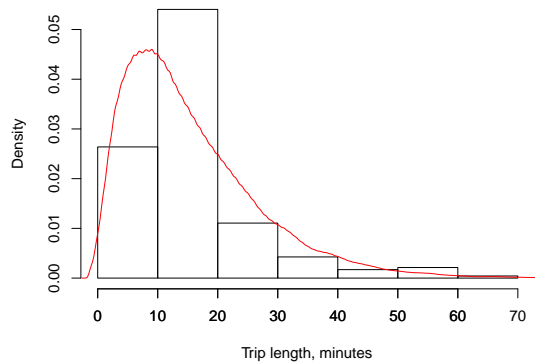}
  \end{center}
  \vspace{-10pt}
  \caption{\emph{The empirical distribution of the trip lengths, shown as the histogram bars, and the theoretical density from the fitted model, shown in red, obtained via Monte Carlo simulation and the subsequent kernel density estimation.}}
  \label{TripLength}
  \vspace{-10pt}
\end{figure}

Fig. \ref{TripLength} represents the histogram of the empirically observed trip lengths along with the theoretical density function of the trip lengths obtained from the fitted model. We use a model with two driving states, as no significant improvement is found beyond this number. Notice that the distribution of the empirically observed trip lengths is adequately captured by the hidden Markov model, although the number of observed trips in the range from 10 to 20 minutes has a higher prevalence than the fitted distribution. In practical applications, more information could be available to model the behavior of the vehicle (e.g. its location and speed), which should facilitate the modeling of the driving patterns.

In the following section, the algorithm for optimally charging the EV is presented. The optimization algorithm makes use of the transition probabilities characterizing the stochastic model for the driving patterns. Thus, the optimization algorithm is designed to handle the stochastic nature of the driving needs.

\section{A Stochastic Dynamic Programming Problem}

The problem of charging an EV can be posed as a conflict between two opposing objectives. The end-user desires to have the vehicle charged and ready for use at his/her discretion, while also minimizing the costs of running the vehicle. Demand for electricity varies over the day and so does the electricity generated from renewable sources. This introduces a varying energy price which can make it beneficial for the end-user to postpone charging his/her vehicle. This means the user is faced with the problem of postponing charging to minimize costs or to charge right away so as to maximize the availability of the vehicle.

The algorithm for optimal charging of the EV is formulated as a stochastic dynamic programming problem. We first define the relevant parameters and variables, and then the state-transition and objective function.

\subsection{Parameters}

\begin{table}[h!]
    \begin{tabular}{ll}
        $u_{\max}$ & maximum rate of charge of the battery (kW)  \\
        $u_{\min}$ & minimum rate of charge of the battery (kW) \\
        $e_{\max}$ & maximum storage level of the battery (kWh)\\
        $e_{\min}$ & minimum storage level of the battery (kWh)\\
        $\lambda_t$      & time varying electricity price (\euro/MWh) \\
        $\phi$     & penalty for violating (unserved) driving needs (\euro/h) \\
        $\eta_c$   & charging efficiency of the battery \\
        $\eta_d$   & discharging efficiency of the battery \\
        $v_i$ 	   & average speed when the vehicle is in use in state $i$ (km/h)\\
        $\mu_i$    & drive efficiency in state $i$ (kWh/km)\\
        $\kappa$   & battery capacity (kWh)\\
        $\omega$   & conversion factor from minutes to hours i.e. $\omega = (60)^{-1}$ (h/min) \\
        $\beta$    & time discount factor\\
\end{tabular}
\vspace{-10pt}
\end{table}
The maximum rate of charge, $u_{\max}$, reflects a power limit on the electric sockets in a residential household or a technical constraint due to thermal limits on the battery (as batteries generate heat when charged). The minimum rate of charge on the battery, $u_{\min}$, reflects that the battery may be limited to only charging i.e. $u_{\min} = 0$ or that discharging the battery is allowed so as to inject power into the grid, i.e. $u_{\min} < 0$. The bounds on the storage limits on the battery, $e_{\max}$ and $e_{\min}$, reflect the storage capacity of the battery. These limits can also be altered to restrict life-cycle degradation of the battery. The penalty $\phi$ is the \emph{inconvenience cost} incurred if the vehicle cannot comply with the driving needs. As seen later, this penalty determines the trade-off between the electricity procurement cost and the availability of the EV to cover a plausible trip. Parameters $\eta_c$ and $\eta_d$ represent the efficiency losses from battery charging and discharging, respectively. The constant $v_i$ is the average speed of the vehicle, when the vehicle is in state $i$, keeping in mind that the modeling framework is general enough to capture multiple different driving states, say urban and rural. The driving efficiency, $\mu_i$, captures the performance of the vehicle in driving state $i$. The constant $\kappa$ is the total energy capacity of the battery. The parameter $\omega$ is used as a conversion factor from hours to minutes, as the model inputs are in hourly values and the model is run in 1-minute time steps.

\subsection{State Variables}

\begin{table}[h!]
    \begin{tabular}{ll}
        $e_t$     & total energy stored in the battery at the beginning of minute $t$ (kWh) \\
        $x_t$     & desired driving state, where $x_t \in \left\{ 1 \ldots N \right\}$ \\
\end{tabular}
\vspace{-10pt}
\end{table}
We assume that variable $x_t$ is exogenously given by the inhomogeneous Markov model described in Section 2. Variable $e_t$ is the energy stored in the battery. We define a state variable at time t as $S_t = (e_t,x_t)$. Notice that, as the driving state is exogenously given, it does not depend on $e_t$ and thus the vehicle is allowed to be in a driving state even though there is no energy on the battery. Logically this is not feasible. Consequently, we refer to $x_t$ as the \emph{desired} driving state, since it can only be reached if there is enough charge on the battery. To cope with this issue, we first define the set, $S_D$, as the collection of states that $x_t$ can take where the vehicle is driving. Then we define the auxiliary variable $x_t^a$ as the actual driving state, i.e.,
\begin{equation}
  x_t^a =
  \left\{
    \begin{array}{ll}
       1 & \textrm{   if } e_t = e_{\min} \wedge x_t \in S_{D} \\
       x_t & \textrm{   else. } \\
    \end{array}
  \right. \label{Definition x_t^a}
\end{equation}
Notice that $x_t$ and $x_t^a$ differ only when there is not enough charge to complete the desired trip. State $1$ denotes one of the parked states. Therefore, according to (\ref{Definition x_t^a}) the vehicle is forced to stop when there is not enough charge on the battery to drive any further. Note that $S_t$ implicitly includes $x^a_t$ as a state, inasmuch as $x^a_t$ is derived from $x_t$ and $e_t$.

\subsection{Decision (Action) Variables}
\begin{table}[h!]
    \begin{tabular}{ll}
        $u_t$     & desired energy charged into (or discharged from) the battery in \\
                  & minute $t$ (kW) \\
\end{tabular}
\vspace{-10pt}
\end{table}
As for the driving state, we define an auxiliary charging variable $u^a_t$, which is the actual energy charged into the battery, since the vehicle is unable to charge when it is in use. The new variable $u^a_t$ is then defined as follows:
\begin{equation}
  u_t^a =
  \left\{
    \begin{array}{ll}
       0 & \textrm{   if } e_t > e_{\min} \wedge x_t \in S_{D} \\
       u_t & \textrm{   else. } \\
    \end{array}
  \right.
\end{equation}
Thus $u_t^a$ is zero when the vehicle is actually driving, and equal to $u_t$ otherwise. Again, if both the state $S_t$ and the desired energy charged $u_t$ are known, the actual energy charged $u_t^a$ follows implicitly from these.

\subsection{State Transition Function}
The driving state variable $x_t$ evolves randomly according to the inhomogeneous Markov model described in Section \ref{Stochastic Model}. The state-transition function for the storage level of the battery can be expressed as:
\begin{equation}
    e_{t+1} = e_t + \left( \eta_c \mathbbm{1}_{ \left\{ u^a_t \geq 0 \right\} } + \frac{1}{\eta_d} \mathbbm{1}_{ \left\{ u^a_t < 0 \right\} }  \right) \omega u^a_t - \sum_{i=1}^N{ \mathbbm{1}_{ \left\{ x^a_t = i \right\} } v_i \mu_i \omega}. \label{Storage dynamics}
\end{equation}
Eq. (\ref{Storage dynamics}) describes the dynamics of the energy stored in the battery. It defines the storage level at time $t+1$, $e_{t+1}$, as the storage level at time $t$, $e_{t}$, plus the net energy charged into the battery and minus the energy that is used to drive the vehicle, which is determined by the random state variable $x_t$. Note that $e_{t+1}$ is written as a function of $e_t$, $x_t^a$ and $u_t^a$. Nevertheless, because $x_t^a$ and $u_t^a$ are functions of $e_t$, $x_t$ and $u_t$, the energy stored in the battery at time $t+1$, $e_{t+1}$, could also be written as functions of these. This, however, would complicate the formulation and for this reason it is omitted here.

\subsection{Constraints}
The desired charging of the battery is limited to being within the bounds for the rate of charge:
\begin{equation}
    u_{\min} \leq u_t \leq u_{\max}. \label{Charge Constraint}
\end{equation}
The storage level on the battery is similarly constrained to being within the storage limits of the battery:
\begin{equation}
    e_{\min} \leq e_t \leq e_{\max}. \label{Energy Constraint}
\end{equation}

\subsection{Objective Function}

The revenue at time period $t$ is given by:
\begin{equation}
    R_t(S_t,u_t) = - \lambda_t \omega u^a_t  - \mathbbm{1}_{ \left\{ x_t \in S_{D} , e_t = e_{\min} \right\} } \omega \phi .
\end{equation}
The first term, $\lambda_t \omega u^a_t$, is the cost incurred from charging the vehicle. The second term, $\mathbbm{1}_{ \left\{ x_t \in S_{D}, e_t = e_{\min} \right\} } \omega \phi$, is the penalty incurred when the user desires to use the vehicle, but he/she cannot do so, because there is not enough energy stored in the battery. Note that this happens precisely when $x_t \neq x_t^a$. Note also that the revenue is equal to the sum of the costs and the penalty with a negative sign.

We introduce now the revenue at the end of the optimization horizon, i.e. at time $T$, when there are no more subsequent decisions to be made:
\begin{equation}
    R_T(S_T) = \eta_d  e_T \frac{1}{T}\sum_{t=1}^{T}{\lambda_t}. \label{Terminal Revenue}
\end{equation}
This equation sets the terminal revenue as the profit that could be made by selling the remaining energy in the battery at the average observed price. One could argue for the use of other terminal conditions: For example, we could replace $\frac{1}{T}\sum_{t=1}^{T}{\lambda_t}$ in (\ref{Terminal Revenue})  with either $\lambda_T$ or $\max_{t}{ \left\{\lambda_t\right\} }$. However, we find it more appropriate to use $\frac{1}{T}\sum_{t=1}^{T}{\lambda_t}$, as this reflects the average economic value of the energy remaining in the battery at $t=T$ if history repeats itself. Besides, $\frac{1}{T}\sum_{t=1}^{T}{\lambda_t}$ constitutes a better prediction of the future electricity price than $\lambda_T$, and $\max_t{\left\{ \lambda_t \right\}}$ would probably lead to an over-estimation of the economic value of the leftover energy, since the battery cannot be fully discharged instantly, even if the maximum electricity price encourages the EV user to do so. A terminal condition is important in obtaining a solution for this problem. However, as explained later, the proposed algorithm is to be applied within a rolling-horizon decision-making process, and as a result, the impact of the terminal condition on the charging pattern is conveniently lessened.

Let $\mathcal{U}_s$ denote the set of feasible decisions according to Eq. (\ref{Definition x_t^a})-(\ref{Energy Constraint}), when the system is in state $s$. Let $\Pi$ denote the set of all feasible policies. A policy, $\pi$, is a collection of decisions $u_t^{\pi}(s) \in \mathcal{U}_s$, spanning the horizon from $t=0$ to $t=T$ and all states $S_t$. Thus for each $t$ and each state $S_t$, $\pi \in \Pi$ will contain the action, $u_t^{\pi}(S_t)$, under the policy $\pi$. For each $\pi \in \Pi$, we can now define the total expected revenue of that policy from time $t$ to $T$ as:
\begin{equation}
    J^{\pi}_t (S_t) = \mathbb{E} \left[  \left. \sum_{\tau=t}^T{ R_{\tau}(S_{\tau}, u^{\pi}_{\tau}(S_{\tau}) ) } \right| S_t \right],
\end{equation}
where $T$ is the optimization horizon. The objective is then to find a policy, $\pi^{\ast}$, that satisfies:
\begin{equation}
    J^{\pi^{\ast}}_t (S_t) = \sup_{\pi \in \Pi}{J^{\pi}_t (S_t)}, \label{Problem Statement}
\end{equation}
for all $0 \leq t \leq T$.

\subsection{Solution Algorithm}

Finding an exact solution to the problem stated in (\ref{Problem Statement}) will be difficult in general due to the randomness and the continuous nature of states and decisions. As the decision at time $t$ depends on the decisions and the values of random variables in previous time periods, the problem grows exponentially as the number of time steps is increased. In order to capture the actual driving patterns and to integrate them into the model in a sensible manner, it is essential that the time resolution is high (1-minute or 5-minute intervals). Due to the fluctuating electricity price and the diurnal variation in the driving pattern, the horizon should be a minimum of one-day ahead. If this is to be accomplished, an exponential growth of the problem is not viable.

Instead we solve the problem by discretizing states and decisions. This yields a discrete stochastic dynamic programming problem that is solved using backward induction and Bellman's principle of optimality. As the driving states are already discrete, the level of energy in the battery and the decision variable $u_t$ remain to be discretized. Suppose that the energy stored in the battery $e_t$ is discretized into $M$ states and there are $N$ driving states. This yields a total of $N \times M$ possible state values for each time step. We define $I_t$ as the index set of possible values that the discretized state variable, $\tilde{S}_t$, can take on. We now define the Bellman equation for the problem stated in Eq. (\ref{Problem Statement}) as
\begin{eqnarray}
    V_t(\tilde{S}_t) &=& \max_{\tilde{u}_t \in \tilde{U}( \tilde{S}_t )} \left\{ R_t(\tilde{S}_t,\tilde{u}_t) + \beta \mathbb{E}_t \left[ V_{t+1}(\tilde{S}_{t+1}) | \tilde{S}_t\right] \right\} \\
    &=& \max_{\tilde{u}_t \in \tilde{U}(\tilde{S}_t)}\left\{ R_t(\tilde{S}_t,\tilde{u}_t) + \beta \sum_{i \in I_{t+1}}{\mathbb{P}_t(\tilde{S}^i_{t+1}|\tilde{S}_t) V_{t+1}(\tilde{S}^i_{t+1})} \right\}, \label{Bellman EQ}
\end{eqnarray}
where $\tilde{S}_t$ is the set of discretized values of $S_t$ and $\tilde{U}(\tilde{S}_t)$ is the set of discretized possible actions in state $\tilde{S}_t$. As the exogenous random variable $X_t$ is defined by a Markov chain, the Bellman equation in Eq. (\ref{Bellman EQ}) represents a Markov decision process, which can be solved using backwards induction, as sketched in Tab. \ref{Pseudocode}.
\begin{table}[h!]
\small
\begin{center}
\begin{tabular}{ll}
  \toprule
  \multicolumn{2}{l}{Backwards Induction Pseudocode}  \\
  \midrule \vspace{-10pt}
  \\
  1: & Initialize: The terminal value is defined as $V_T(\tilde{S}_T) = R_T(\tilde{S}_T)$ given by (\ref{Terminal Revenue})\\
  2: & \textbf{for} $t = T-1$ to 0 \textbf{do}: \\
  3: & \quad    \textbf{for} $s \in I_t$ \textbf{do}:  \\
  4: & \qquad        $\tilde{u}^{\ast}_{s,t} = \textrm{arg} \max_{\tilde{u}_{t} \in \tilde{U}(s)} \left\{ R_t(s,\tilde{u}_{t}) + \beta \sum_{i \in I_{t+1}}{\mathbb{P}_t(i|s) V_{t+1}(i)} \right\}$ \\
  5: & \qquad $V_t(s) = R_t(s,\tilde{u}^{\ast}_{s,t}) + \beta \sum_{i \in I_{t+1}}{\mathbb{P}_t(i|s) V_{t+1}(i)} $ \\
  6: & \quad  \textbf{end} $s$ \textbf{for} \\
  7: & \textbf{end} $t$ \textbf{for}  \\
  \vspace{-10pt} \\
  \bottomrule
\end{tabular}
    \caption{\label{Pseudocode} Pseudocode using backwards induction to obtain the optimal policy $\tilde{\pi}^{\ast}$ as the collection of $\tilde{u}^{\ast}_{s,t}$.}
\end{center}
\vspace{-20pt}
\end{table}

Using the algorithm in Tab. \ref{Pseudocode}, we find the optimal discretized policy $\tilde{\pi}^{\ast}$ as the collection of $\tilde{u}^{\ast}_{s,t}$ indicating when and how much to charge, depending on both the time $t$ and the pair $s$ of driving and battery states.

It is advisable to run the algorithm over a long horizon, say two days, to incorporate the diurnal variation in the driving pattern and in the energy price. In addition, the longer the horizon covered by the optimization process, the smaller the influence of the terminal condition. Indeed, we propose to re-run the algorithm in Tab. \ref{Pseudocode} for every time step, with a horizon that is extended accordingly, following a rolling-window process.

\section{Results and Discussion}
In this section the model has been implemented and run with realistic parameter values. The data at our disposal only includes two states, \emph{driving} and \emph{not driving}. We consider a Markov model with three states; one time-varying \emph{not driving} state and two time-invariant \emph{driving} states. Model results are compared with those obtained from ``rule-of-thumb'' policies to assess the economic performance of the proposed decision-support tool. We first present an in-sample study with the model fitted to the training set, which serves to illustrate its main features. We then carry out an out-of-sample study to evaluate the performance of the model on the test set. For simplicity, we assume that the vehicle is plugged into the electricity grid when not driving.

\subsection{Model Characteristics}
We consider an EV with a battery capacity $\kappa = 24 \; \rm{ kWh}$ and an average consumption of $\mu_i = 0.20 \; \rm{ kWh/km}$. The entire battery capacity is assumed to be available for use, i.e. $e_{\max} = 24 \; \rm{ kWh}$ and $e_{\min} = 0 \; \rm{ kWh}$. We also consider that the vehicle is mainly to be employed in an urban driving cycle with an average speed of $v_i = 40 \; \rm{ km/h}$, including stopping for red lights and congestion. This yields a range of 120 km on one charge and a drive time of 3 hours. Regarding the charging, we assume a maximum charging capacity of $u_{\max} = 4 \; \rm{ kW}$. In the base case the vehicle is not allowed to discharge power back into the grid, i.e. $u_{\min} = 0 \; \rm{ kW}$. This case is subsequently extended to allow for discharging via a V2G scheme with $u_{\min} = -4 \; \rm{ kW}$. The charging efficiency parameters are $\eta^c = \eta^d = 0.9$. On the basis of these characteristics, the vehicle resembles the Nissan Leaf, which is one of the top selling EVs in the world (as of January 2013).

We consider an optimization horizon covering 48 hours in advance to incorporate the diurnal variation in the energy price as well as in the driving pattern. Furthermore, as already explained in Section 3.7, the relatively long horizon is used to decrease the influence of the terminal condition on the optimal charging scheme, which is gradually obtained from the rolling-window process. The time resolution in the model is in minutes, which yields a total of 2880 time steps. A 1-minute time resolution is chosen to adequately model the use of the vehicle. As we consider a horizon of 48 hours, the discount factor is set to $\beta = 1$. The state variable for the energy charged on the battery is discretized into 360 different states. Likewise, the state variable for the use of the vehicle has three states (one \emph{not driving} and two \emph{driving} states). Therefore, the model relies on $3 \times 360$ different states for each time step. In the base case, where only charging is allowed, the vehicle charges at either full rate or not at all, thus the decision variable $u_t$ can only take on two different values. The optimal solution is found in less than a minute on a personal computer with a 2.70 GHz processor and 8.0 GB RAM, which is satisfactory. The model can be straightforwardly modified to work with 5-minute or 10-minute time resolution with a view to further decreasing the solution time. Also, the discretization of the energy charged on the battery can be coarser. This may be useful if the model is extended with more driving states, or the model has to be implemented with less computing capacity, or if the optimization horizon has to be extended. However, as the model run-time is quite small, such efforts have not been pursued. We notice that the model is parameterized in terms of the penalty, $\phi$, incurred when the vehicle does not have enough energy in the battery to complete the desired trip. This can also be seen as a risk-aversion parameter, where the risk of not completing a trip is weighed against minimizing the costs of driving.

With regard to the electricity price, we use the Nordpool DK1 spot-price historical series. We consider that the EV charging controller receives a 48-hour price forecast from, for example, a distribution system operator (DSO). The Nordpool spot price is determined each day in blocks of 24 hourly values and is made public at noon the previous day. Therefore we assume that the risk associated with the volatility of the electricity price is handled by the DSO or some other intermediary, but not by the end-consumer. Besides, the 48-hour price forecast may be updated, if appropriate, every time (every minute) the model is re-run as part of the rolling-window process.

\subsection{In-Sample Study}

We use next the training data set defined in Section 2.2 to fit the stochastic model for driving patterns. Then we simulate plausible driving scenarios based on this model and evaluate the performance of the proposed decision-support tool for optimal charging on these scenarios. Therefore, the analysis carried out here is in-sample, i.e. it assumes that the fitted stochastic model for driving patterns perfectly captures the actual nature of the use of the vehicle. The purpose of this study is then to illustrate the main features of the proposed decision-making tool. Firstly, we analyze schemes where only charging is permitted. Then we consider V2G schemes~\cite{kempton1997electric, Kempton2005268}, where the vehicle is permitted to supply power from the battery to the grid. We use electricity prices from 00:00 on the 25-01-2012 to 00:00 on the 29-01-2012.

\subsubsection{Charging-Only Schemes }

\begin{figure}[h]
  \begin{center}
    \includegraphics[width=0.8\textwidth]{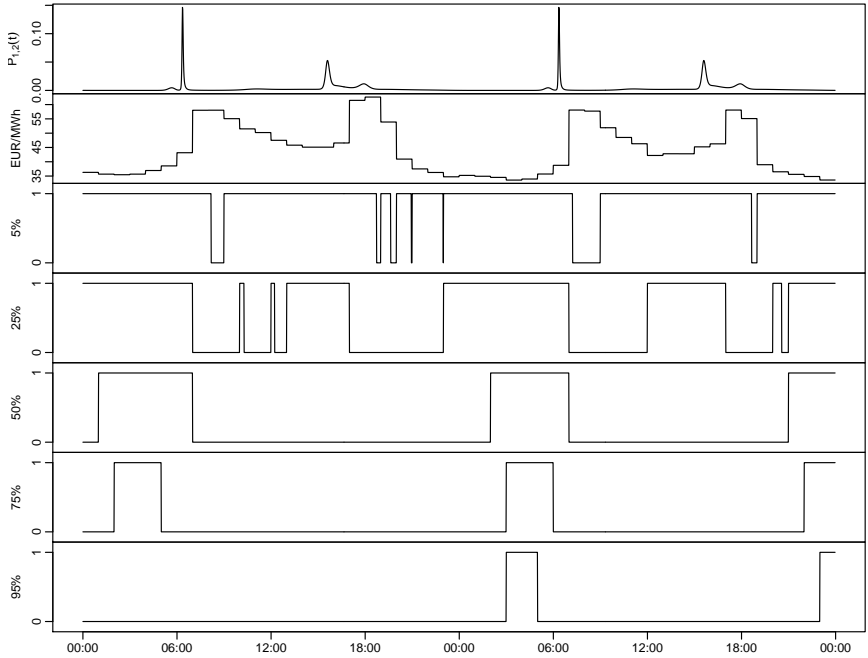}
  \end{center}
  \vspace{-20pt}
  \caption{ \emph{From top to bottom: time-varying probabilities of starting a trip, electricity price, and $\tilde{\pi}^{\ast}$ for different levels of charge on the battery and penalty $\phi = 10 $.}}
  \vspace{-10pt}
  \label{Charging}
\end{figure}

Fig. \ref{Charging} shows the estimated time-varying probability of starting a trip, the electricity price, and selected values for the optimal policy $\tilde{\pi}^{\ast}$, which defines the appropriate charging action to be undertaken given the state $\tilde{S}_t$ and the time $t$. The optimal policy may take values in the set $\left\{1 , 0 \right\}$ for charging and not charging, respectively. In Fig. \ref{Charging} the battery state is indicated on the vertical axis for different levels of charge, expressed as a percentage of the battery capacity $e_{\max}$. The time is indicated on the horizontal axis. Additionally, note that Fig. \ref{Charging} only shows the charging decisions for when the vehicle is not driving, as we assume that it is not possible to stop a trip and recharge, unless the battery is fully depleted. It is important to stress that Fig. \ref{Charging} shows a single run of the optimization algorithm. The difference between this snap-shot of the algorithm and the rolling window process will be illustrated subsequently.

It can be observed from this figure that if the energy level of the battery is $5\%$, the optimal decision is to always charge, except in those time periods when the probability of driving is low and the electricity price is particularly high. In contrast, if the energy level of the battery is 50\%, the vehicle is only charged when the energy price is comparatively low. This charging policy becomes more extreme as the level of charge approaches 100\%. Indeed, if the energy level of the battery is equal to 95\% of $e_{\max}$, the EV is only charged in those time periods where the energy price is expected to reach its lowest values.

\begin{figure}[h]
  \begin{center}
    \includegraphics[width=0.8\textwidth]{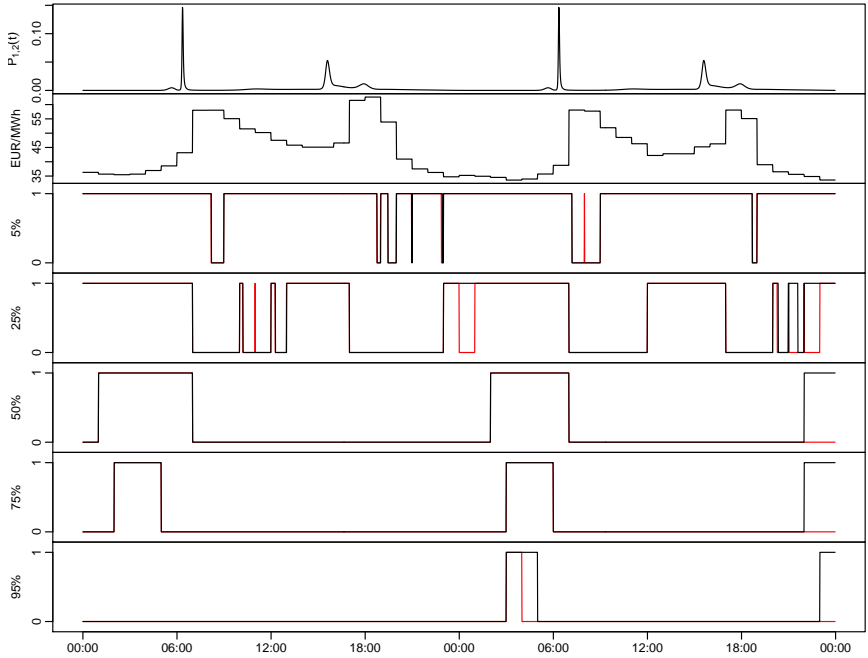}
  \end{center}
  \vspace{-20pt}
  \caption{\emph{Analogous to Fig. \ref{Charging}. A rolling-horizon policy has been implemented and the resulting policy is shown in red along with the fixed-window policy in black in the bottom five graphs.}}
  \label{RollingHorizon}
  \vspace{-10pt}
\end{figure}

In reality, the proposed charging algorithm is to be used following a rolling-horizon process, which allows for updating the energy price forecast and reducing the effect of the terminal condition on the optimal policy, as highlighted next. In Fig. \ref{RollingHorizon} the results yielded by the algorithm when implemented over a fixed two-day horizon are compared to those obtained considering a two-day rolling-horizon. In the rolling-horizon optimization, the model is rerun every hour and the optimal policy updated accordingly. The rolling-horizon is kept fixed to two days in advance, and consequently we use energy prices from 00:00 on the 25-01-2012 to 00:00 on the 29-01-2012, that is, four days in total. From Fig. \ref{RollingHorizon} we see that there are only slight deviations between the rolling horizon and the fixed horizon procedures within the first day. On the second day, however, we begin to see deviations that go beyond a single spike. As the time approaches time $T$, more discrepancies are observed between the two models, indicating that the terminal condition has an impact on the optimal charging policy. However, this impact is mostly confined to the last time periods of the 48-hour optimization horizon, and therefore implementation of the proposed algorithm in a rolling-horizon fashion can reduce, if not completely eliminate, this effect.

\begin{figure}[h]
  \begin{center}
    \includegraphics[width=0.8\textwidth]{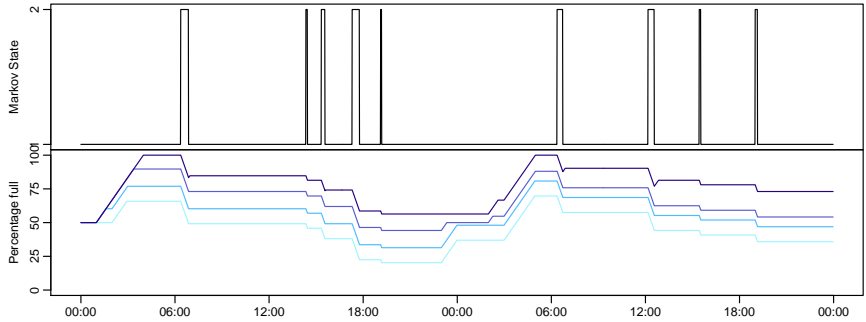}
  \end{center}
  \vspace{-20pt}
  \caption{\emph{Top: realization of the driving pattern. Bottom: the corresponding charge on the battery in percent for different penalty values when implementing $\tilde{\pi}^{\ast}$. The lightest blue line refers to the lowest penalty and the penalty increases with the darker shades of blue.}}
  \label{ParametersSimulation}
  \vspace{-10pt}
\end{figure}

To illustrate the actual implementation of an optimal policy $\tilde{\pi}^{\ast}$, we run the following simulation process. Firstly, $\tilde{\pi}^{\ast}$  is computed and defined for different values of the penalty $\phi$. Secondly, a plausible realization of the driving pattern is simulated and the corresponding time evolution of the level of charge on the battery is determined according to this realization and the optimal policy $\tilde{\pi}^{\ast}$. The results of such a simulation are shown in Fig. \ref{ParametersSimulation}. It can be seen that as the penalty increases, the level of charge on the battery is correspondingly higher, conditional on the same realization of the driving pattern.

Another promising aspect of EVs is the possibility of supplying power into the grid at times of high demand. This is investigated in the following section.

\subsubsection{Vehicle-to-Grid Schemes}

Allowing for the vehicle to supply power from the battery into the grid has the potential to help mitigate the effects of peak power demand. This operation mode is usually referred to as a Vehicle-to-Grid (V2G) scheme. A V2G scheme is investigated here from the perspective of a single vehicle.

\begin{figure}[h!]
  \begin{center}
    \includegraphics[width=0.8\textwidth]{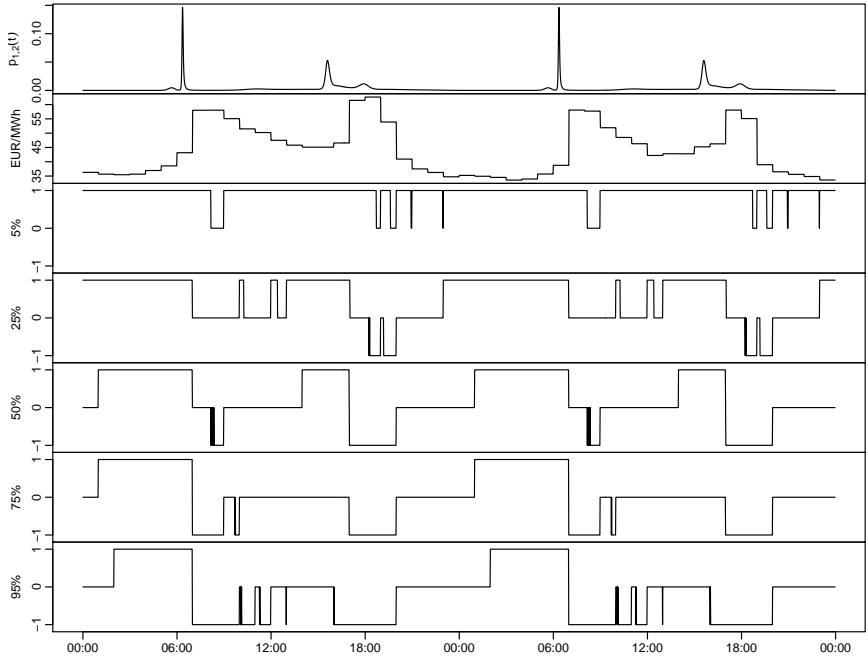}
  \end{center}
  \vspace{-20pt}
  \caption{\emph{Top to bottom: time-varying probabilities of starting a trip, electricity price, and $\tilde{\pi}_{V2G}^{\ast}$ for different levels of charge on the battery with penalty $\phi = 10$. In this case, the discharge of the EV battery back into the grid is considered and therefore, the policy can take on negative values. This figure is similar to Fig. \ref{RollingHorizon}, where only charging is allowed.}}
  \label{Discharging}
  \vspace{-10pt}
\end{figure}

Implementation of the V2G scheme is by setting $u_{\min}$ to $-4 \rm{ kW}$ and keeping all other parameter values unchanged. The optimal policy, $\tilde{\pi}^{\ast}_{V2G}$, obtained by implementing a V2G scheme for a penalty value of $\phi = 10$, is shown in Fig. \ref{Discharging}, which is similar to Fig. \ref{RollingHorizon}, except that the optimal policy may take values in the set $\left\{1 , 0 ,-1 \right\}$ for charging, not charging, and discharging, respectively. Observe that when the energy level in the battery is low, the optimal policy, $\tilde{\pi}^{\ast}_{V2G}$, basically consists of charging the whole time, except in those periods with electricity prices at their peak and a low probability of driving. As the energy level in the battery increases, the policy changes to supplying power into the grid at the price peaks and to charging at the price valleys. The proposed V2G algorithm thus weighs the costs associated with running out of charge on the battery against the gains from delaying charging to when the energy price is low and the gains from supplying power into the grid when the energy price is high. Also in Fig. \ref{Discharging}, notice that the optimal policy shows some ``spikes'' where the optimal decision is changed from charging to not charging for a short time. This is linked to the fact that the electricity price traded on Nordpool is in hourly time resolution, and therefore the price changes only every hour and the corresponding price change may be large. As the vehicle decides the appropriate action for every minute, it is able to exploit this in its charging strategy.

\subsection{Out-of-Sample Study}

We now evaluate the model performance on the test data set defined in Section 2.2. Therefore, we provide results from testing the optimal charging policies on the actual utilization of the vehicle in the second half of the data period. This study is thus performed out of sample. We use electricity prices from 00:00 on 01-01-2011 to 00:00 on 08-03-2011.

\begin{figure}[h]
  \begin{center}
    \includegraphics[width=0.8\textwidth]{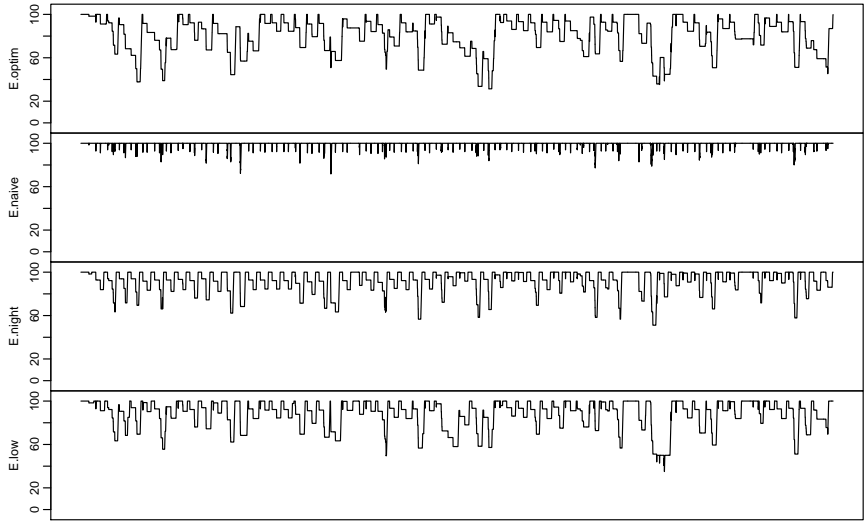}
  \end{center}
  \vspace{-30pt}
  \caption{ \emph{State of charge in \% for different charging policies, from top to bottom: Optimized charging with $\phi = 10$ (with only charging allowed), naive charging, night charging and low-price charging. }}
  \vspace{-10pt}
  \label{RealDrivingDataCharging}
\end{figure}

Fig. \ref{RealDrivingDataCharging} shows the state of charge of the battery for different charging policies. Note that this figure is analogous to the bottom plot of Fig. \ref{ParametersSimulation}, except that different charging policies are considered and the time shown is 65 days. In particular, the top graph in Fig. \ref{RealDrivingDataCharging} is obtained using the proposed decision-support tool to find the optimal policy for charging (V2G operation mode is not permitted). The lower three graphs represent different ``rule of thumb'' policies. We refer to the first one as ``naive charging'', according to which the vehicle is charged immediately upon being parked. The second to last is called ``night charging'', and entails charging the vehicle at night between 10 pm and 6 am or if the charge on the battery is below 50 \%. The last one is ``low price charging'', under which the vehicle is charged if the electricity price is in the lowest 20\%-quantile of the price distribution for the next 24 hours or if the charge on the battery is below 50 \%. From Fig. \ref{RealDrivingDataCharging}, notice that none of the strategies empties the battery at any time.

\begin{figure}[h]
  \begin{center}
    \includegraphics[width=0.8\textwidth]{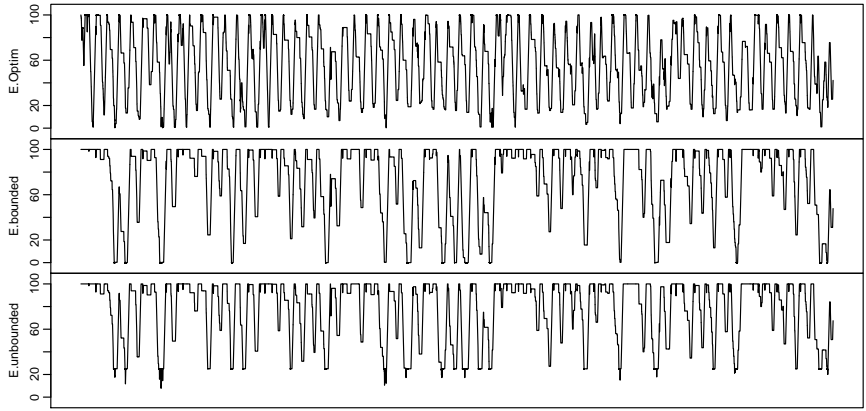}
  \end{center}
  \vspace{-30pt}
  \caption{ \emph{State of charge in \% for different charging policies, from top to bottom: Optimized charging with $\phi = 10$ (with V2G charging allowed), ``rule of thumb'' V2G unbounded, ``rule of thumb'' V2G bounded. }}
  \vspace{-10pt}
  \label{fig:RealDrivingDataV2G}
\end{figure}

Fig. \ref{fig:RealDrivingDataV2G} is similar to Fig. \ref{RealDrivingDataCharging}, except that in this case the V2G operation mode is allowed. The optimal charging policy is compared with two other ``rule of thumb'' policies in which V2G is permitted. More specifically, in the ``rule of thumb'' V2G charging schemes, the vehicle is charged if the electricity price is in the lowest 30\%-quantile of the price distribution for the next 24 hours and discharged if the electricity price is above the 90\%-quantile. We consider an unbounded policy, with no restriction on the lower limit of charge on the battery, and a bounded policy, where this lower limit is set to 25\% of the battery capacity. This lower limit specifies when the vehicle is available for supplying power into the grid.

Comparing Fig. \ref{RealDrivingDataCharging} and Fig. \ref{fig:RealDrivingDataV2G}, it is seen that the total battery capacity is exploited when using a V2G scheme. This has the side effect of the battery being depleted at some times during the day. Consequently, the vehicle is not able to cover the user's driving needs, should the user desire to drive.

Tab. \ref{tab:OutOfSample} compares the costs and availability of the vehicle under the different charging policies. We consider the average daily cost of running the vehicle in \euro \ and the number of events where the vehicle is not able to cover the user driving needs counted over the 65 days that the test data set spans.

\begin{table}[h!]
\small
\begin{center}
    \begin{tabular}{rrrrr}
    \toprule
          & \multicolumn{2}{c}{Charging only} & \multicolumn{2}{c}{V2G permitted}\\
    Penalty, $\phi$      & Cost, \euro & Events & Cost, \euro & Events \\ \midrule
    2                    & 0.170 & 0  & -0.097  & 12  \\
    5                    & 0.174 & 0  & -0.084  & 8   \\
    10                   & 0.177 & 0  & -0.061  & 3   \\
    100                  & 0.181 & 0  & -0.047  & 0   \\
    1000                 & 0.188 & 0  & -0.019  & 0   \\ \midrule
    Naive                & 0.323 & 0  & -      &  -   \\
    Night                & 0.284 & 0  & -      &  -   \\
    Low Price            & 0.210 & 0  & -      &  -   \\
    V2G unbounded        & -    & - & 0.071   &  11   \\
    V2G bounded          & -    & - & 0.133   &  0   \\ \bottomrule
\end{tabular}
\end{center}
\vspace{-10pt}
  \caption{\emph{Average daily costs in \euro \ and number of events where there is enough charge on battery to service user driving needs.}}
\label{tab:OutOfSample}
\end{table}

Let us consider first the strategies in Tab. \ref{tab:OutOfSample} under which only charging is allowed, we see that there are no observed events of not having enough charge on the battery to complete a trip. Also, we notice, as expected, that the optimal charging strategies have lower costs than the ``rule of thumb'' policies. The low-price charging strategy is indeed the ``rule of thumb'' policy that approximates closest to the optimal policy in terms of costs and availability. It yields, however, an average daily cost which is around 12-24\% higher than that obtained from implementing the proposed decision-support tool.

As for the charging policies that include V2G operation mode, it becomes apparent that caution should be exercised to prevent the vehicle from being fully discharged when the end-user desires to drive. In this line, notice that increasing the penalty reduces the number of observed events of not having enough charge on the battery to cover a desired trip. Introducing a V2G charging scheme allows for substantially reducing the cost associated with driving as opposed to charging-only schemes, and may even result in negative average costs. Observe that the optimal charging policy developed in this paper clearly outperforms the ``rule of thumb'' V2G schemes. In the unbounded case, charging costs are substantially reduced, but multiple out-of-battery events are recorded. Imposing a lower bound on the discharging solves this problem, but at the expense of considerably increasing the running cost of the vehicle, to such an extent that it nearly doubles.

The difference in performance between the optimal charging strategy and the ``rule of thumb'' policies can be expected to become larger for electric vehicles covering higher distances or with lower battery capacity.

Lastly, we would like conclude this section by pointing out that, in general, the spot price is not the price observed by the end-user. Indeed, the end-user faces a price that includes taxes and other costs on top of the spot electricity price. As an example, consider a country like Denmark, where the average price of electricity paid by the end-user, including taxes and fees, is around  \euro 300 /MWh, which is 5-10 times the average spot price \citep{http://www.energy.eu}. In the current Danish power system, fees and taxes are imposed on the amount of electricity consumed by the end-user, not on its total cost. This does not encourage the end-user to switch to a smart consumption of energy based on variable prices. In fact, if the taxes and fees were implemented as a function of the total energy cost, the savings from switching to a smart charging policy in Denmark could be multiplied by a factor of between five and ten.

\section{Conclusion}

This paper proposes an algorithm to optimally charge an electric vehicle based on stochastic dynamic programming. The algorithm is built on an inhomogeneous (hidden) Markov chain model that characterizes the stochastic use of the vehicle. The algorithm determines the optimal charging policy depending on the use of the vehicle, the risk aversion of the end-user, and the electricity price. The costs associated with running the vehicle are decreased significantly when the charging strategy is determined by the proposed optimization model, with little or no inconvenience to the end-user. These costs can be reduced even further if the vehicle is permitted to supply power into the grid. Indeed, findings show the possibility of making a net profit from running the vehicle. The proposed stochastic dynamic programming model is versatile and can easily be adapted to any specific vehicle, thus providing a customized charging policy.

A possible extension would be to apply the proposed model to data with more Markov states, which could be used to investigate the benefits of installing more public charging stations as opposed to home charging, or to capture different driving states such as ``urban'', ``rural'', or ``highway''. In addition, the model could be enhanced to consider transition probabilities that are estimated adaptively in time. An adaptive approach would capture structural changes in the driving behavior, such as variations over the year or a change in use that could follow, for example, from the householder buying an additional vehicle. Moreover, adaptivity is relevant for applying the model in practice.

Further research could be also directed at modeling a fleet of vehicles by using a mixed-effects model. The optimization scheme could be applied individually to each vehicle and the total population load could be evaluated. This would highlight if and how EVs could be used to mitigate an increase in peak electricity demand when switching from combustion-based vehicles to EVs. Other investigations could focus on the relationship between EVs and renewable energy sources and how EVs could be used to move the excess production to time periods of high demand, possibly making renewables more economically competitive.

\section*{Acknowledgments}
DSF (Det Strategiske Forskningsr\aa d) is to be acknowledged for partly
funding the work of Emil B. Iversen, Juan M. Morales and Henrik Madsen through the Ensymora project (no. 10-093904/DSF). Furthermore, Juan M. Morales and Henrik Madsen are partly funded by the iPower platform project, supported by DSF (Det Strategiske Forskningsr\aa d) and RTI (R\aa det for Teknologi og Innovation), which are hereby acknowledged. Finally, we thank DTU Transport for providing the data used in this research.


\section*{References}




\end{document}